\documentclass[12pt,draftcls,onecolumn]{IEEEtran}
\usepackage{amsthm}
\usepackage{amsmath,amssymb}			

\newtheorem{thm}{Theorem}[section]
\newtheorem{theorem}[thm]{Theorem}
\newtheorem{Theorem}[thm]{Theorem}

\newtheorem{corollary}[thm]{Corollary}

\newtheorem{lemma}[thm]{Lemma}

\theoremstyle{definition}

\newtheorem{example}{Example}[section]

\theoremstyle{remark}

\newcommand{\RR}{\mathbb{R}}                                  

\begin{document}

\title{On Delay-independent Stability of a class of Nonlinear Positive Time-delay Systems}
\author{Vahid S. Bokharaie \thanks{Corresponding Author: Vahid S. Bokharaie is with Mathematics Applications Consortium for Science and Industry (MACSI), Univeristy of Limerick, Ireland. Email: vahid.bokharaie@ul.ie}, Oliver Mason \thanks{Oliver Mason is with Hamilton Institute, National University of Ireland Maynooth, Ireland. Email: oliver.mason@nuim.ie}}
\maketitle


\begin{abstract}                          
We present a condition for delay-independent stability of a class of nonlinear positive systems.  This result applies to systems that are not necessarily monotone and extends recent work on cooperative nonlinear systems.  
\end{abstract}

\section{Introduction} \label{sec:intro}
In Ecology, Biology, Economics and other application domains, the variables of interest are typically constrained to be nonnegative.  Motivated by this very simple observation, numerous researchers have worked on the theory of so-called positive systems in the recent and not so recent past \cite{FarRin00}.  Much of the more recent work in this direction has focussed on extending the elegant theory of positive linear time-invariant systems to more general and realistic settings \cite{RKW, ELENA, NGOC, NGOC2}.  In particular, authors have considered time-varying, switched and different classes of nonlinear positive systems \cite{Kha99, KNG99, Ric03}.  

Of the particular properties of positive LTI systems, that of delay-independent stability is among the most striking.  Formally, if the positive system 
\begin{equation}
\label{eq:zerodellin}
\dot{x}(t) = (A+B)x(t) 
\end{equation}
has a globally asymptotically stable (GAS) equilibrium at the origin, then so does the delayed system
\begin{equation}
\label{eq:delaylin} \dot{x}(t) = A x(t) + B x(t - \tau) 
\end{equation}
for every $\tau \geq 0$.  

Recently, this result has been extended to nonlinear cooperative systems in \cite{MV, BMV}.  In a similar spirit, the robustness of other classes of positive systems with respect to delay has been studied in \cite{NGOC}.  In this brief note, our purpose is to build on the work of \cite{BMV} and derive a condition for delay-independent stability of a class of nonlinear positive systems that are not necessarily cooperative.  To this end, we present necessary preliminary results in Section \ref{sec:pre}; the main result of the note is then developed in Section \ref{sec:main} along with a simple illustrative example.  

\section{Preliminaries} 
\label{sec:pre}
Throughout we denote the positive orthant of $\mathbb{R}^n$ by $\mathbb{R}^n_+ := \{x \in \RR^n: x_i \geq
0, 1 \leq i \leq n\}.$   For $x \in \RR^n$ and $i = 1, \ldots , n$, $x_i$ denotes the $i^{th}$
coordinate of $x$.  For vectors $x, y \in \RR^n$, we write: $x \geq y$ if
$x_i \geq y_i$ for $1 \leq i \leq n$; $x
> y$ if $x \geq y$ and $x \neq y$; $x \gg y$ if $x_i > y_i, 1 \leq i \leq
n$. 

Let $D$ be an open neighbourhood of the nonnegative orthant $\mathbb{R}^n_+$ and let $f:D \rightarrow \mathbb{R}^n$, $g:D \rightarrow \mathbb{R}^n$ be $C^1$ mappings.  We are concerned with the delayed system 
\begin{equation}\label{eq:del1}
\dot{x}(t) = f(x(t)) + g(x(t - \tau)), \;\;\; \tau \geq 0.
\end{equation}
For a continuous function, $\phi \in C([-\tau, 0] , \mathbb{R}^n)$ we write $F(\phi) = f(\phi(0)) + g(\phi(-\tau))$ for the functional associated with (\ref{eq:del1}). Also, $x(t, \phi)$ denotes the solution of (\ref{eq:del1}) corresponding to the initial condition $\phi \in C([-\tau, 0], \mathbb{R}^n)$.  Throughout the note, we assume that the origin is an equilibrium of (\ref{eq:del1}) so that $(f+g)(0) = 0$.  For background on the theory of delayed systems and functional differential equations, see \cite{Hale}.

The system (\ref{eq:del1}) is \emph{positive} if for any initial condition $\phi \in\mathcal{C}:= C([-\tau, 0], \mathbb{R}^n_+)$, the solution $x(t, \phi) \in \mathbb{R}^n_+$ for all $t \geq 0$ for which it is defined.  It follows from Theorem 5.2.1 of \cite{SMITH} that (\ref{eq:del1}) will be positive if 
\begin{equation}
\label{eq:POS}
F_i(\phi) \geq 0 \mbox{ for all } \phi \in \mathcal{C} \mbox{ with } \phi_i(0) = 0.
\end{equation}
This condition will be satisfied if: 
\begin{itemize}
\item[(P1)] $g(x) \geq 0$ for all $x \in \mathbb{R}^n_+$;
\item[(P2)] $f_i(x) \geq 0$ for all $x \in \mathbb{R}^n_+$ with $x_i = 0$.  
\end{itemize}
For the remainder of the note, unless clearly stated otherwise, we assume that $f$ and $g$ satisfy (P1), (P2).   As all the systems considered are positive, we shall work with the state spaces $\mathcal{C}$ (for delayed systems) and $\mathbb{R}^n_+$ (for undelayed systems) throughout.

 Moreover we shall assume that $f$ and $g$ are subhomogeneous of degree $\alpha > 0$ meaning that $f(\lambda x) \leq \lambda^\alpha f(x)$ for all $\lambda \geq 1$, $x \in \mathbb{R}^n_+$. The class of subhomogeneous vector fields given above includes concave vector fields \cite{Kra01}. Furthermore, it includes vector fields which are homogeneous with respect to the standard dilation map

\textit{KKM Lemma}

Later in the paper, we shall need the so-called Knaster, Kuratowski,
Mazurkiewicz (KKM) Lemma \cite{KKM2}.  We denote the simplex whose vertices are the standard basis vectors $e_1, \ldots, e_n$ of $\mathbb{R}^n$ by $\Delta_n$.  Given a set of indices $1 \leq i_1 < i_1 < \cdots < i_p \leq n$, the
simplex $S(e_{i_1}, \ldots , e_{i_p})$ with the vertices $e_{i_1}, \ldots , e_{i_p}$ is a \textit{face} of $S(e_1, \ldots , e_n)$.  We shall need
the following open version of the KKM Lemma.  

\begin{Theorem}[KKM Lemma]
\label{thm:KKM}
Let $F_1, \ldots , F_r$ be
(relatively) open subsets of $\Delta_n$.  If 
\begin{equation}\nonumber
 S(e_{i_1}, \ldots , e_{i_p}) \subset F_{i_1} \cup \cdots \cup F_{i_p}
\end{equation}
holds for all faces $S(e_{i_1}, \ldots , e_{i_p})$, $1 \leq p \leq n$, $1 \leq i_0 < i_1 < \cdots
< i_p \leq n$, then 
\begin{equation}\nonumber
F_1 \cap \cdots \cap F_{n} \neq \emptyset. 
\end{equation}
\end{Theorem}

\textit{Quasimonotone Conditions and Monotonicity}

Our main result in the next section, which applies to positive systems that are not necessarily monotone, makes extensive use of the properties of monotone systems and conditions for monotonicity.  We now recall some relevant definitions and results.

A function $h: \mathbb{R}^n_+ \rightarrow \mathbb{R}^n$ satisfies the quasimonotone condition (on $\mathbb{R}^n_+$) if $h_i(x) \leq h_i(y)$ for any $x, y \in \mathbb{R}^n_+$ satisfying $x \leq y$, $x_i = y_i$.  It is well known \cite{SMITH} that this implies that the associated system $\dot{x}(t) = h(x(t))$ is monotone, meaning that $x_0 \leq x_1$ implies $x(t,  x_0) \leq x(t, x_1)$ for all $t$ for which both solutions are defined.  

For a positive delayed system (\ref{eq:del1}), the quasimonotone condition requires that $F_i(\phi) \leq F_i(\psi)$ for any $\phi, \psi$ in $\mathcal{C}$ with $\phi \leq \psi$ and $\phi_i(0) = \psi_i(0)$.  Again this is a sufficient condition for (\ref{eq:del1}) to be monotone \cite{SMITH}, meaning that $\phi \leq \psi$ implies $x(t, \phi) \leq x(t, \psi)$ for all $t$ for which both solutions are defined.  

Given a vector $v \in \mathbb{R}^n_+$, we use $\hat{v}$ to denote the function in $\mathcal{C}$ with $\hat{v}(t) = v$ for $t \in [-\tau, 0]$. 

\section{Main Result}
\label{sec:main}
In this section, we develop the main result of this note.  First we present some preliminary technical lemmas.  

\begin{lemma}
\label{lem:KKM} Let $h: \mathbb{R}^n_+ \rightarrow \mathbb{R}^n$ be continuous.  Further, assume that $h_i(x) \geq 0$ for any $x \in \mathbb{R}^n_+$ with $x_i = 0$ and that there exists no $x \neq 0$ in $\mathbb{R}^n_+$ with $h(x) \geq 0$.  Then there exists some $v \in \Delta_n$, with $h(v) \ll 0$. 
\end{lemma}
\textbf{Proof:} We apply the KKM lemma, suitably adapting the arguments given in \cite{BMV, DRW}.  For $1 \leq i \leq n$, let 
$$F_i := \{ x \in \Delta_n : h_i(x) < 0 \}.$$
As $h$ is continuous, $N_i$ is relatively open in $\Delta_n$ for $1 \leq i \leq n$.  Consider indices $1 \leq i_1 < i_2 \cdots  < i_p \leq n$ and let $x \in S(e_{i_1}, \ldots , e_{i_p})$ be given.  It follows from the assumptions on $h$, that $h_i(x) < 0$ for some index $i$.  However $x_j = 0$ and hence $h_j(x) \geq 0$ for $j \notin \{i_1, \ldots , i_p\}$ so it follows that $h_{i_j}(x) < 0$ for some $j \in \{1, \ldots, p\}$.  This implies that 
$$ S(e_{i_1}, \ldots , e_{i_p}) \subset F_{i_1} \cup \cdots \cup F_{i_p}.$$  The KKM Lemma now implies the result.

The vector $v$ whose existence is established in the previous result will play a key role in the stability analysis presented in Theorem \ref{thm:Main}.  The approach to stability taken here is reminiscent of the so-called MO condition used to analyse Wazewski systems (see \cite{Alex1, Mart} and the references therein).  

The next result is the main result of this note and provides a sufficient condition for the origin to be a GAS equilibrium of (\ref{eq:del1}) for all $\tau \geq 0$.  

\begin{theorem}
\label{thm:Main} Consider the positive system (\ref{eq:del1}) and suppose that $f$ and $g$ are subhomogeneous of degree $\alpha > 0$.  Assume that for every $w \in \mathbb{R}^n_+ \setminus \{0\}$, there is some index $i$ such that 
\begin{equation}
\label{eq:Cond1} 
\textrm{sup}\{g_i(y) : 0 \leq y \leq w \} < - \textrm{sup} \{f_i(x) : 0 \leq x \leq w, x_i = w_i \}.
\end{equation}
Then the origin is a GAS equilibrium of (\ref{eq:del1}) for every $\tau \geq 0$. 
\end{theorem}
\textbf{Proof:}  
We shall adapt the techniques used to prove Proposition 5.2.3 of \cite{SMITH}; note that we cannot directly apply this result as we are not explicitly assuming that (\ref{eq:del1}) possesses an invariant order interval.  We associate with (\ref{eq:del1}) a positive time-delay system 
\begin{equation}
\label{eq:CompSys} \dot{x}(t) = \bar{F}(x_t) 
\end{equation}
with the following properties:
\begin{itemize}
\item the system (\ref{eq:CompSys}) is order-preserving;
\item the trajectories of (\ref{eq:del1}) are dominated by those of (\ref{eq:CompSys}) for every $\tau \geq 0$;
\item the origin is a GAS equilibrium of (\ref{eq:CompSys}) for every $\tau \geq 0$.
\end{itemize}
Taken together, these points will yield the desired result.  Following \cite{SMITH} for any $\phi \in \mathcal{C}$ define
\begin{equation}
\label{eq:comp} \bar{F}_i(\phi) = \textrm{sup} \{ F_i(\psi) : 0 \leq \psi \leq \phi, \psi_i(0) = \phi_i(0) \}, \quad \; 1 \leq i\leq n.
\end{equation}
Remembering that $F(\psi) = f(\psi(0))+ g(\psi(-\tau))$, it is readily seen that
\begin{eqnarray*}
\bar{F}_i(\phi) &=& \textrm{sup} \{ f_i(x) + g_i(y) : 0 \leq x \leq \phi(0), x_i = \phi_i(0), 0 \leq y \leq \phi(-\tau) \}.
\end{eqnarray*} 
As both $f_i$ and $g_i$ are continuous and the set $\{(x, y) \in \mathbb{R}^{2n} : 0 \leq x \leq \phi(0), x_i = \phi_i(0), 0 \leq y \leq \phi(-\tau)\}$ is clearly compact, it follows that $\bar{F}$ is well-defined and in fact the supremum is a maximum.  Using this last observation, and noting that both $f$ and $g$ are $C^1$ and hence Lipschitz on any compact set, it is straightforward to show directly that $\bar{F}$ is Lipschitz on compact subsets of $\mathcal{C}$. 

We next verify that $\bar{F}$ satisfies the positivity requirement (\ref{eq:POS}).  To this end, let $\phi \in \mathcal{C}$ with $\phi_i(0) = 0$ be given.  As $F$ satisfies (\ref{eq:POS}) by assumption, it follows that for any $\psi \in \mathcal{C}$ with $\psi_i(0) = \phi_i(0) = 0$,$\psi \leq \phi$, we must have $F_i(\psi) \geq 0$.  Hence $\bar{F}_i(\psi) \geq 0$ and $\bar{F}$ satisfies (\ref{eq:POS}).  

As $\bar{F}$ is Lipschitz on compact subsets of $\mathcal{C}$ and satisfies (\ref{eq:POS}), for any $\phi \in \mathcal{C}$ there exists a unique solution $x(t, \phi)$ to (\ref{eq:CompSys}), which satisfies $x(t, \phi) \in \mathcal{C}$ for all $t$ in its maximal interval of existence.  

It is immediate from the definition of $\bar{F}$ that for any $\phi \in \mathcal{C}$, $F(\phi) \leq \bar{F}(\phi)$.  The argument from \cite{SMITH} to establish that $\bar{F}$ satisfies the quasimonotone condition applies directly.  In the interests of completeness, we outline it now.  Let $\phi \leq \psi$ in $\mathcal{C}$ and suppose $\phi_i(0) = \psi_i(0)$ for some $i$.  Then the supremum defining $\bar{F}_i(\phi)$ is taken over a subset of that defining $\bar{F}_i(\psi)$ and hence $\bar{F}_i(\phi) \leq \bar{F}_i(\psi)$.  Thus $\bar{F}$ satisfies the quasimonotone condition and hence the system (\ref{eq:CompSys}) is monotone. It is clear that (\ref{eq:CompSys}) has an equilibrium at the origin.   

Next define $h:\mathbb{R}^n_+ \rightarrow \mathbb{R}$ by $h(v) = \bar{F}(\hat{v})$.  We claim that there exists no non-zero $w \in \mathbb{R}^n_+$ with $h(w) \geq 0$.  To see this, let $w \neq 0$ in $\mathbb{R}^n_+$ be given and note that for each $i$
$$h_i(w) = \textrm{max}\{ f_i(x) + g_i(y) : 0 \leq x, y \leq w, x_i = w_i\}.$$
By assumption, there is some $i$ such that (\ref{eq:Cond1}) holds.  For this $i$, we must have $h_i(w) < 0$.  Thus, we have shown that for every non-zero $w \in \mathbb{R}^n_+$, there exists some index $i$ with $h_i(w) < 0$.  This immediately implies that the origin is the unique equilibrium of (\ref{eq:CompSys}) in $\mathbb{R}^n_+$.  

For $v \in \mathbb{R}^n_+$, $1 \leq i \leq n$:
$$h_i(v) = \textrm{sup}\{f_i(x) + g_i(y) : 0 \leq x \leq v, x_i = v_i , 0 \leq y \leq v \}.$$ 
As $f$ and $g$ are both subhomogeneous of degree $\alpha$ by assumption, it is readily seen that $h$ is also subhomogeneous.  

Clearly, for any $x \in \mathbb{R}^n_+$ with $x_i = 0$, we must have $h_i(x) \geq 0$.  We have also shown above that there exists no $w \neq 0$ in $\mathbb{R}^n_+$ with $h(w) \geq 0$.  It now follows from Lemma \ref{lem:KKM} that there must exist some vector $v \gg 0$ with $h(v) = \bar{F}(\hat{v}) \ll 0$.  Corollary 5.2.2 of \cite{SMITH} now implies that the solution $x(t, \hat{v})$ of (\ref{eq:CompSys}) converges to an equilibrium of (\ref{eq:CompSys}) as $t \rightarrow \infty$.  However, as noted above, 0 is the only such equilibrium.  Thus $x(t, \hat{v}) \rightarrow 0$ as $t \rightarrow \infty$.  

Now for any $\phi \in \mathcal{C}$, there exists some $\lambda > 1$ such that $\phi \leq \lambda v$.  Moreover, as $h$ is subhomogeneous, $\bar{F}(\lambda \hat{v}) \leq \lambda^{\alpha} \bar{F}(\hat{v}) \ll 0$, the solution $x(t, \lambda \hat{v}) \rightarrow 0$ as $t \rightarrow \infty$.  Hence as (\ref{eq:CompSys}) is monotone and positive, it follows that $x(t, \phi) \rightarrow 0$ as $t \rightarrow \infty$.  This completes the proof.  

We next present a simple example to illustrate the above result.
\begin{example}
\label{eg:Pos1} Let $f:\mathbb{R}^2_+ \rightarrow \mathbb{R}^2$, $g:\mathbb{R}^2_+ \rightarrow \mathbb{R}^2$ be given by 
$$ f(x_1, x_2) = \left( \begin{array}{c} 
							x_1(1-e^{x_1+x_2}) \\
							-x_2
						\end{array}\right), 
\;\;
g(x_1, x_2) = \left( \begin{array}{c}
						x_1x_2 \\
						\frac{x_2}{1+x_2}
					\end{array} \right).
$$
Then it is easy to see that $f$ and $g$ are subhomogeneous of degree 2 and that $F(\phi) = f(\phi(0)) + g(\phi(-\tau))$ satisfies the positivity requirement (\ref{eq:POS}).  Note also that $\frac{\partial f_1}{\partial x_2}$ is not non-negative for all $x_1, x_2 \geq 0$ so the system is not monotone.  

Now let $w = (w_1, w_2)^T \in \mathbb{R}^n_+ \setminus \{0\}$ be given with $w_2 > 0$.  It is easy to see that 
$\textrm{sup}\{g_2(x): 0 \leq x \leq w\} = \frac{w_2}{1+w_2}$.  Also, $\textrm{sup}\{f_2(y): 0 \leq y \leq w, y_2 = w_2\} = -w_2$.  Hence as $w_2 > 0$, 
$$\textrm{sup}\{g_2(x): 0 \leq x \leq w\} < - \textrm{sup}\{f_2(y): 0 \leq y \leq w, y_2 = w_2\}.$$
Next suppose $w_2 = 0$.  Then, $w_1 > 0$ as $w \neq 0$ by assumption, and 
$\textrm{sup}\{g_1(x): 0 \leq x \leq w\} = 0$ and $\textrm{sup}\{f_1(y): 0\leq y \leq w, y_1 = w_1\} = w_1(1-e^{w_1}) < 0$.  Hence, in this case, 
$$\textrm{sup}\{g_1(x): 0 \leq x \leq w\} < - \textrm{sup}\{f_1(y): 0 \leq y \leq w, y_1 = w_1\}.$$
It now follows from Theorem \ref{thm:Main} that the system (\ref{eq:del1}) has an asymptotically stable equilibrium at the origin for every value of $\tau > 0$.
\end{example}

As a final point for this section, we note the following simple corollary for monotone time-delay systems.  The system (\ref{eq:del1}) will be monotone if $f$ is cooperative  and $g$ is non-decreasing (see Chapter 5 of \cite{SMITH}). 

\begin{corollary}
\label{cor:monotone} Consider the system (\ref{eq:del1}) and assume that $f$ is cooperative and $g$ is non-decreasing.  Assume that the system
\begin{equation}
\label{eq:undel1} \dot{x}(t) = (f+g)(x(t))
\end{equation}
has a globally asymptotically stable equilibrium at the origin.  Then the system (\ref{eq:del1}) also has a GAS equilibrium at the origin for all $\tau > 0$. 
\end{corollary}
\textbf{Proof:}  As the origin is GAS, there can be no non-zero $w \in \mathbb{R}^n_+$ with $(f+g)(w) \geq 0$.  Thus for every such $w$, there is some index $i$ with $f_i(w) + g_i(w) < 0$.  As $g$ is non-decreasing, we must have $g_i(y) \leq g_i(w)$ for all $y \leq w$ in $\mathbb{R}^n_+$.  Also, as $f$ is cooperative, it follows that $f_i(x) \leq f_i(w)$ for all $x \in \mathbb{R}^n_+$ with $x \leq w$, $x_i = w_i$.  Thus 
\begin{eqnarray*}
 \textrm{sup}\{g_i(x): x \leq w\}  & \leq & g_i(w) \\
 & <& -f_i(w) \\
 &\leq & - \textrm{sup} \{ f_i(y):  y \leq w, y_i = w_i \}.
\end{eqnarray*}
It now follows immediately from Theorem \ref{thm:Main} that (\ref{eq:del1}) has a GAS equilibrium at the origin for all $\tau \geq 0$.

\section{Conclusions}
\label{sec:con}
We have extended some recent work on delay-independent stability for positive systems.  Specifically, in Theorem \ref{thm:Main}, a sufficient condition for a class of nonlinear positive systems to be stable independent of delay is given.  The systems covered by the result are not necessarily monotone.  In fact, the corresponding result for monotone delay systems follows as a simple corollary. 
\section*{Acknowledgements}
This work was supported by the Irish Higher Educational Authority (HEA) PRTLI 4 Network Mathematics Grant and Science Foundation Ireland award 09/SRC/E1780.



\begin{thebibliography}{99.}
\bibitem{FarRin00}
Lorenzo Farina and Sergio Rinaldi.
\newblock {\em Positive Linear Systems. Theory and Applications}.
\newblock Pure and Applied Mathematics. John Wiley \& Sons, Inc., New York, NY,
  USA, 2000.

\bibitem{RKW}
Bj\"{o}rn~S. R\"{u}ffer, Christopher~M. Kellett, and Steven~R. Weller.
\newblock Connection between cooperative positive systems and integral
  input-to-state stability of large-scale systems.
\newblock {\em Automatica}, 46(6):1019--1027, 2010.

\bibitem{ELENA}
Ettore Fornasini and Maria-Elena Valcher.
\newblock Linear copositive {L}yapunov functions for continuous-time positive
  switched systems.
\newblock {\em IEEE Transactions on Automatic Control}, 55(8):1933--1937, 2010.


\bibitem{NGOC}
P.~H.~Anh Ngoc, S.~Murakami, T.~Naito, J.~Son Shin, and Y.~Nagabuchi.
\newblock On positive linear {V}olterra-{S}tieltjes differential systems.
\newblock {\em Integral Equations and Operator Theory}, 64:325--255, 2009.

\bibitem{NGOC2}
P.~H.~Anh Ngoc.
\newblock Stability of positive differential systems with delay.
\newblock {\em IEEE Transactions on Automatic Control}, 58(1):203--209, 2013.

\bibitem{Kha99}
Kharitonov, V. L.
\newblock {\em Robust stability analysis of time delay systems: A survey}.
\newblock {\em Annual Reviews in Control}, 23:185-196, 1999.

\bibitem{KNG99}
Kolmanovskii, V. B. and Niculescu, S. I. and Gu, K.
\newblock {\em Delay effects on stability: A survey}.
\newblock {\em Proceedings of the 38th IEEE Conference on Decision and Control}, 2:1993-1998, 1999.

\bibitem{Ric03}
Richard, J. P.
\newblock {\em Time-delay systems: an overview of some recent advances and open problems}.
\newblock {\em Automatica}, 39(10):1667--1694, 2003.

\bibitem{MV}
Oliver Mason and Mark Verwoerd.
\newblock Observations on the stability properties of cooperative system.
\newblock {\em Systems and Control Letters}, 58:461--467, 2009.


\bibitem{BMV}
Vahid~S. Bokharaie, Oliver Mason, and Mark Verwoerd.
\newblock D-stability and delay-independent stability of homogeneous
  cooperative systems.
\newblock {\em IEEE {T}ransactions on {A}utomatic {C}ontrol},
  55(12):1996--2001, 2010.

\bibitem{Hale}
J.~K Hale and S.~M. Verduyn~Lunel.
\newblock {\em Introduction to Functional Differential Equations}.
\newblock Springer-Verlag, 1993.

\bibitem{Kra01}
Ulrich Kraus.
\newblock {\em Concave Perron-Frobenius theory and applications}.
\newblock {\em Nonlinear Analysis-Theory Methods and Applications}, 47(3):1457--1466, 2001.

\bibitem{KKM2}
Marc Lassonde.
\newblock Sur le principe {K}{K}{M}.
\newblock {\em Comptes rendus de l'Acad\'{e}mie des sciences}, 310(7):573--576,
  1990.


\bibitem{SMITH}
Hal~A. Smith.
\newblock {\em Monotone Dynamical Systems. An introduction to the Theory of
  Competitive and Cooperative Systems}.
\newblock American Mathematical Society, Providence, Rhode Island, USA, 1995.


\bibitem{DRW}
Sergey Dashkovskiy, Bj\"{o}rn~S. R\"{u}ffer, and Fabian~R. Wirth.
\newblock An {I}{S}{S} small-gain theorem for general networks.
\newblock {\em Mathematics of Control, Signals and Systems}, 19:93--122, 2007.


\bibitem{Alex1}
A. Yu. Aleksandrov, A. V. Platonov.
\newblock {\em On Stability and Dissipativity of Some Classes of Complex Systems}.
\newblock	{\em Automation and Remote Control}, 70(8):1265--1280, 2009.

\bibitem{Mart}
A. A. Martynyuk.
\newblock {\em Asymptotic Stability Criterion for Nonlinear Monotonic Systems and its Applications}.
\newblock {\em International Applied Mechanics}, 47(5):3--67, 2011.


\end{thebibliography}
\end{document}